\documentclass[12pt]{article}

\usepackage{amsmath,amssymb,amsmath,amsthm,xspace,a4wide,amscd}

\theoremstyle{plain}

\newcommand{\dps} {\displaystyle}

\newcommand{\rr} {\rightarrow}
\newcommand{\C}{\mathbb{C}}

\newcommand{\x} {\mbox{$x\kern-0.90em-$}}

\newcommand{\la}{\longrightarrow}

\newcommand{\N}{{\mathbb{N}}}

\newcommand{\R}{\mathbb{R}}
\newcommand{\HE}{\cal{H}}
\newcommand{\ope}{{\cal{O}}}

\newcommand{\sa} {S^A(\R^n)}
\newcommand{\CBA} {CB^{\infty}(\R^n,A)}
\newcommand{\CtwoRA}{CB^\infty(\R^{2n},A)}

\newcommand{\op} {{\cal{B}}^{*}(E)}
\newcommand{\dc}{d\!\!\! /}

\newcommand{\epsi}{\varepsilon}

\newcommand{\ai}{\hfill \rule [0.01in]{0.07in}{0.07in}} 

\newtheorem{prop}{{\bf Proposition}}[section]

\newtheorem{lema}{{\bf Lemma}}[section]

\newtheorem{teo}{{\bf Theorem}}[section]

\newtheorem{nota}{{\bf Note}}[section]

\newtheorem{ob}{{\bf Remark}}[section]

\newtheorem{corol}{{\bf Corollary}}[section]

\newcommand{\pf}{\noindent{\sc Proof.}\,\,}

\numberwithin{equation}{section}

\begin{document}

\date{}












%


\title{{\Large BOUNDEDNESS OF PSEUDODIFFERENTIAL OPERATORS 
OF $C^*$-ALGEBRA-VALUED SYMBOL.}\thanks{During the completion of this work, the author received finantial support from CAPES, Brazil,
through its "conv\^enio PRO-DOC". 2000 AMS Classification: Primary 47G30.}}

\author{Marcela I. Merklen}

\maketitle

\thispagestyle{empty}


\abstract{

\noindent
Let us consider the set $\sa$ of rapidly decreasing functions $G:\R^n \rightarrow A$, 
where $A$ is a separable $C^*-$algebra. We prove a version
of the Calder\'on-Vaillancourt theorem for pseudodifferential operators
acting on $\sa$ whose symbol is $A$-valued. 
Given a skew-symmetric matrix, $J$, we prove that a pseudodifferential operator
that commutes with $G(x+JD)$, $G\in\sa$, is of the form $F(x-JD)$, 
for $F$ a $C^\infty-$function with bounded derivatives of all orders.

\bigskip

\section{\bf Introduction.}

\vspace{0,5cm}

Throughout this work, $A$ denotes a separable $C^*-$algebra
and $\sa$ denotes the $A-$valued Schwartz space of smooth and rapidly decreasing
functions on $\R^n$. On $\sa$ we define the $A$-valued inner product

\[<f,g> = \dps\int f(x)^*g(x)dx,\] whose associated norm we denote by $\|\cdot\|_2$,
$\|f\|_2 = \|<f,f>\|^{\frac{1}{2}}$.

\medskip The completion of $\sa$ with this norm is a Hilbert $A-$module that
we denote by $E$. The set of all adjointable (and therefore bounded) operators on $E$ we denote by $\op$. Let $\CBA$ denote the set of $C^\infty$-functions with bounded derivatives
of all orders.

\bigskip

\medskip
In section 2, 
we see a generalization of the Calder\'on-Vaillancourt Theorem, \cite{c-v}, for a pseudodifferential 
operator, $a(x,D)$,  whose symbol, $a$, is in $CB^\infty(\R^{2n},A)$. 

Following the steps of Seiler's proof in \cite{sei-1} (in fact going back to Hwang, \cite{hw}),
we can see that $a(x,D)$ is bounded on
$E$. 
Note that Seiler's result needs that one works on a Hilbert space where we have
that the Fourier transform is unitary in the usual sense. One of the advantages of working
with this norm, $\|\cdot\|_2$, is that the Fourier transform becomes an "unitary" operator on
the completion of $\sa$; {\em i. e.} it is a Hilbert-module adjointable operator on $E$, whose
inverse is equal to its adjoint. This norm $\|\cdot\|_2$ allows for a proof of the Calder\'on-Vaillancourt Theorem for 
pseudodifferential operators whose symbols are $A$-valued functions.

We also prove that these operators are adjointable.
So, we have $a(x,D)\in \op$, $a\in CB^\infty(\R^{2n},A)$.


\bigskip
In \cite{ri}, Rieffel defines a deformed product in $\CBA$, depending on an
anti-symmetric matrix, $J$, by \[F\times_JG(x) = 
\dps\int (2\pi)^{-n}e^{iu\cdot v}F(x+Ju)G(x+v)dudv.\]
It is not difficult to see that the left-regular representation of $\CBA$ defines 
pseudodifferential operators on $\op$; in fact, for $F\in\CBA$, 
$L_F(\varphi) =F\times_J\varphi$, $\varphi\in\sa$, is the
pseudodifferential operator of symbol $F(x-J\xi)$.

\medskip
Rieffel proves that $L_F, \,\, F\in \CBA $ is a continuous operator on $E$,
(\cite{ri}, Corollary 4.7) and that $L_F$ is adjointable on $E$
(\cite{ri}, Proposition 4.2).

\medskip
The Heisenberg group acts on $\op$ by conjugation in the
following way.

\medskip
Given $V\in\op$, 
\[ (z,\zeta,t)\la E_{z,\zeta,t}^{-1} V E_{z,\zeta,t}, \hspace{1cm} 
(z,\zeta,t)\in\R^{2n}\times\R,\] where 
\[ E_{z,\zeta,t}f(x) = e^{it}e^{i\zeta x}f(x-z), \hspace{1cm} f\in\sa.\]

It is easy to see that $V_{z,\zeta} = E_{z,\zeta,t}^{-1}VE_{z,\zeta,t}$
does not depend on $t\in\R$.

\medskip
We say that $V$ is {\em Heisenberg-smooth} if the map $(z,\zeta)\la V_{z,\zeta}$ is 
$C^\infty$, and, if $z\la V_z$ is $C^\infty$, where $V_z = V_{z,0}$,
we say that $V$
is {\em translation-smooth}.

When we are dealing with the scalar case, $A=\C$, we have the remarkable
characterization of Heisenberg-smooth operators in $\op$ given by H. O. Cordes, \cite{co-2}:
these are the psudodifferential operators whose symbols are in $CB^\infty(\R^{2n})$.

\bigskip

In section 3 we prove that if a skew-symmetric, $n\times n$, matrix is given
and if the $C^*-$algebra $A$ is such that a suitably stated generalization of Cordes' 
characterization can be proved, then any Heisenberg-smooth operator $T\in\op$,
which commutes with every psudodifferential operator with symbol $G(x+J\xi)$, for some
$G\in\CBA$, is also a psudodiferential operator with symbol $F(x-J\xi)$, for some $F\in\CBA$.
This is a rephrasing of a conjecture stated by Rieffel for an arbitrary $A$ at the end of
Chapter 4 of \cite{ri} (the operators $G(x+JD)$ are those obtained from the right regular 
representation for his deformed product on $\CBA$). That Cordes' characterization implies Rieffel's conjecture has already been proved for the scalar case, \cite{m-m}. The
Schwartz kernel argument used in \cite{m-m} has to be avoided here, in the more general case.




\section{\bf $a(x,D)\in\op$.}

Let us consider a pseudodifferential operator on $E$ such that, if $\varphi\in\sa$,

\[a(x,D)\varphi(x) = \int e^{i(x-y)\xi}a(x,\xi)\varphi(y)\dc y\dc\xi,\] for
$a\in \CtwoRA$, where $\dc y = (2\pi)^{-\frac{n}{2}}dy$. 
As in the scalar case, we can see that $a(x,D)\varphi(x)$ is well defined 
for each $x\in\R^n$, if $\varphi\in\sa$.

\medskip
An example of such an operator is given by Rieffel:

Given a function $F\in\CBA$,\[L_F\varphi(x) = 
\int e^{i(x-y)\xi}F(x-J\xi)\varphi(y)\dc y\dc\xi.\] 
The integrals considered here are {\em oscilatory integrals} (\cite{ri}, Chapter 1).

\vspace{1cm} Let us see next the fundamental ideas of a generalization of the 
Calder\'on-Vaillancourt Theorem for operators on $E$. 

 \bigskip
First, let us see that $a(x,D)(\varphi)\in E$, for $\varphi\in\sa$.

\medskip
Considering $L^2(\R^n,A) $ as the set of all functions $f:\R^n\rightarrow A$
such that $\int \|f(x)\|^2dx<\infty$, with the ``almost everywhere" equivalence
relation, where we consider the norm $\|\cdot\|_{L_2}$ (defined in the
usual way), we can prove that $a(x,D)(\varphi)\in L^2(\R^n,A)$, as follows:

\medskip
Using integration by parts and the equation 
\begin{equation}\label{E:i+x}
 (i+x)^\alpha e^{ixy} = (i+D_y)^\alpha e^{ixy}, \hspace{0.8cm} \alpha = (1,...,1), 
\hspace{0.5cm} \text{where} \,\,\, D_y = -i\partial_y,\end{equation}

\noindent we obtain

\[a(x,D)(\varphi)(x) = \dps\int e^{ix\xi}a(x,\xi)\hat{\varphi}(\xi)\dc\xi = \]
\[ = (i+x)^{-\alpha}\dps\int\left[(i+D_\xi)^{\alpha} e^{ix\xi}
\right]a(x,\xi)\hat{\varphi}(\xi)\dc\xi =\]
\[= (i+x)^{-\alpha}\dps\int e^{ix\xi}\left[(i-D_\xi)^\alpha a(x,\xi)\hat{\varphi}(\xi)\right]\dc\xi,\]
where $\hat{\varphi}(\xi) = \dps\int e^{-i\xi y}\varphi(y)\dc y$. Since
the last integral is bounded, $a(x,D)(\varphi)\in L^2(\R^n,A)$.

\medskip 
On the other hand, since $\|\cdot\|_2 \leq \|\cdot\|_{L2}$, we can prove that $L^2(\R^n,A)\subseteq E$, so that $a(x,D)(\phi)\in E$.

\begin{teo}Let $a\in CB^\infty(\R^{2n},A)$. Then, $a(x,D)$ is a bounded 
operator on $E$. In fact, \,$\|a(x,D)\|\leq l\pi(a)$, for $l\in\R^+$ independent of $a$, 
and $\pi(a)=
\dps\sup\{\|\partial^\beta_x\partial^\gamma_\xi a\|_\infty \,\,; \\ \beta,
\gamma\leq \alpha = (1,1,\cdots,1)\}$.
\end{teo}

\pf
To begin with, let us consider the case when $a$ has compact support. Denoting 
$a(x,D)$ by $T$,  for  $\varphi,\psi\in\sa,$ we look at 
$<\psi,T\varphi>$, which equals $<\hat{\psi},\widehat{T\varphi}>$.
(Here we are
dealing with the Fourier transform in $L^2(\R^n,A)$, which is ``unitary" on $E$, in the sense
that $<f,g> = <\hat{f},\hat{g}>, \hspace{1cm} f, g\in L^2(\R^n,A)\subseteq E.)$

Since 
\[\widehat{T\varphi}(\eta) = \int e^{-i\eta x}T\varphi(x)\dc x = 
\int e^{-i\eta x}e^{i(x-y)\xi} a(x,\xi)\varphi(y)\dc y\dc\xi\dc x,\]

\medskip
we have

\[ <\hat{\psi},\widehat{T\varphi}> = 
\int e^{-i\eta x}e^{i(x-y)\xi}\hat{\psi}^*(\eta)a(x,\xi)\varphi(y)\dc y\dc \xi\dc x d\eta.\]

\medskip
Using integration by parts and the equation (\ref{E:i+x}), we have

\[  <\hat{\psi},\widehat{T\varphi}> = \int e^{-i\eta x}(i+x-y)^{-\alpha}\hat{\psi}^*(\eta)
(i+x-y)^\alpha e^{i(x-y)\xi}  a(x,\xi)\varphi(y)\dc y\dc \xi\dc x d\eta = \]

\[ = \int e^{-ix\eta}(i+x-y)^{-\alpha}\hat{\psi}^*(\eta) e^{ix\xi} 
e^{-iy\xi}\left[(i-D_\xi)^\alpha a(x,\xi)\right]\varphi(y) \dc y\dc \xi\dc x d\eta= \]

\[ = \int e^{-iy\xi}(i+x-y)^{-\alpha}(i+\xi-\eta)^{-\alpha}\hat{\psi}^*(\eta) e^{ix(\xi-\eta)}
(i+\xi-\eta)^\alpha\left[(i-D_\xi)^\alpha a(x,\xi)\right]\varphi(y) \dc y\dc \xi\dc x d\eta= \]
\[ = \int e^{ix\xi}e^{-ix\eta}(i+x-y)^{-\alpha}\hat{\psi}^*(\eta) e^{-iy\xi}
(i+\xi-\eta)^{-\alpha}\left[(i-D_x)^\alpha(i-D_\xi)^\alpha a(x,\xi)\right]\varphi(y)
\dc \eta\dc y\dc x d\xi.\]


\bigskip
Let us consider \[h(z) = (i-z)^{-\alpha} \hspace{1cm} \alpha=(1,\cdots,1)\] and
\[h_z(y) = h(y - z)  \hspace{1cm} y,z\in\R^{n}.\] 


Thus, we have

\[<\hat{\psi},\widehat{T\varphi}> = \int e^{ix\xi}f(x,\xi)\left[(i-D_x)^\alpha(i-D_\xi)^\alpha
a(x,\xi)\right]g(x,\xi)\dc x d\xi,\] with 

\[ f(x,\xi) = \int e^{-ix\eta}h_\xi(\eta)\hat{\psi}(\eta)^*\dc\eta\] and

\[g(x,\xi) = \int e^{-iy\xi}h_x(y)\varphi(y)\dc y.\] So, we can write (by
abuse of notation),
\[ \|<\hat{\psi},\widehat{T\varphi}>\| = (2\pi)^{-\frac{n}{2}}\|<e^{-ix\xi}f^*(x,\xi),
\left[(i-D_x)^\alpha(i-D_\xi)^\alpha a(x,\xi)\right]g(x,\xi)>\|.\] If $c(x,\xi) = 
(i-D_x)^\alpha(i-D_\xi)^\alpha a(x,\xi)$, there exists $d_1\in\R^+$, not depending on $a$,
such that $\dps\sup_{(x,\xi)\in\R^{2n}}\|c(x,\xi)\| < d_1\pi(a).$

In Prop. \ref{exist-c} below, we prove that there exists 
$d_2\in\R^+$, not depending on $\varphi$ or $g$, such that $\|g\|_2 \leq d_2\|\varphi\|_2$,
so that we have 
\[ \|cg\|_2 \leq d_1\pi(a)d_2\|\varphi\|_2.\] 


\medskip 

In a similar way as in Proposition \ref{exist-c},
we get that
\[\|f^*\|_2 \leq d_2\|\psi\|_2,\] then, for $k=d_1d_2^2(2\pi)^{-\frac{n}{2}}$, 
we have, for all $\varphi,\psi\in\sa$,
\[\|<\psi,a(x,D)\varphi>\| \leq k\pi(a)\|\varphi\|_2\|\psi\|_2.\] 

\medskip
As for the general case, we consider the function $a_\epsi\in CB^\infty(\R^{2n},A)$
given by $a_\epsi(x,\xi) = \phi(\epsi x,\epsi \xi)a(x,\xi)$, for $0<\epsi\leq 1$, and
where $\phi\in C^\infty_c(\R^{2n})$ is such that $\phi \equiv 1$ close to zero.
As we just have seen, we have that $\|<\psi,a_\epsi(x,D)\varphi>\|
\leq k\pi(a_\epsi)\|\varphi\|_2\|\psi\|_2$. 
Just doing some computations, we get
that there is $m\in\R^+$, not depending on $a$ or $\epsi$, such that 
$\pi(a_\epsi)\leq m\pi(a)$. Besides, it is not difficult
to see that we have $\lim_{\epsi\rightarrow 0}<\psi,a_\epsi(x,D)\varphi> \,=\,
<\psi,a(x,D)\varphi>$. 
Actually,
\[\|<\psi,a(x,D)\varphi>\|\underset{\epsi\rightarrow 0}\longleftarrow \|<\psi,a_\epsi(x,D)\varphi>\|\leq
km\pi(a)\|\varphi\|_2\|\psi\|_2. \] 

\medskip
\noindent Considering, now, $a(x,D)\varphi$ in place of $\psi$, since
$a(x,D)\varphi\in L_2(\R^n,A)$, we have that
\[\|<a(x,D)\varphi,a(x,D)\varphi>\| \leq l\pi(a)\|a(x,D)\varphi\|_2\|\varphi\|_2,
\hspace{1cm} l=km, \, \forall\varphi\in\sa\]
as before. So, there is $l\in\R^+$, not depending on $a$, such that $\|a(x,D)\|\leq l\pi(a)$.

\hfill$\ai$

\begin{prop}\label{exist-c} Given $\varphi\in\sa$, let $g(x,\xi)=\dps\int e^{-iy\xi}
(i+x-y)^{-\alpha}\varphi(y)\dc y$, then there exists $d\in\R^+$, not depending neither
on $g$ nor on $\varphi$, such that $\|g\|_2\leq c\|\varphi\|_2$.\end{prop}

\pf Let $h_x$ be as before, and put $g(x,\xi)=\dps\int e^{-iy\xi}h_x(y)\varphi(y)\dc  y = 
\widehat{h_x\varphi}(\xi)$. Then

\[\int g(x,\xi)^*g(x,\xi)dxd\xi = \int<\widehat{h_x\varphi},\widehat{h_x\varphi}>dx =
\int<h_x\varphi,h_x\varphi>dx \] 
\[=\int \overline{h(x)}h(x)dx\int\varphi(\xi)^*\varphi(\xi)d\xi .\] If $d = \left(\dps\int |h(x)|^2dx\right)^{\frac{1}{2}}$, we have $\|g\|_2\leq d\|\varphi\|_2. \,\,\, \ai$

\begin{nota}{\rm If $a\in CB^\infty(\R^{2n},A)$, we denote by 
${\cal{O}}(a)$  the pseudodifferential
operator whose symbol is $a$.}\end{nota}

\begin{prop}\label{exist-p}  
There exists $p\in \CtwoRA$ such that 
\begin{equation}\label{E:calO}
<\ope(a)\varphi,\psi> = <\varphi,\ope(p)\psi> \hspace{1cm}  \hspace{1cm}
\forall \varphi,\psi\in\sa.\end{equation}                      \end{prop}

\pf  \,\, First \, we \, prove \, that \hspace{0,1cm}
$p(y,\xi) = \dps\int e^{-ix\eta}a(y-z,\xi-\eta)^*\dc z\dc\eta $ \hspace{0,1cm} belongs \, to  \, $\CtwoRA$.
As for this, we  use strongly the definition of oscilatory integrals given in\cite{ri},
where we consider for a while the corresponding Fr\'echet space $\CtwoRA$.

Then, applying proposition 1.6 of \cite{ri}, we can begin 
working with $a$ of compact support, for which we can work as Cordes in 
chapter 1 section 4 of \cite{co-1}.

To obtain the general case, we apply the Dominated Convergence Theorem.
Please, see details at proposition 4.6 of \cite{ire}. \,\, $\ai$

%

\begin{ob}{\rm   The application $\ope:\CtwoRA\la\op$, given by $a \mapsto a(x,D)$, 
is well defined and it is easy to see that it is injective. 
  }\end{ob}

\begin{ob}\label{scalar}{\rm   As in the scalar case, \cite{co-1}, chapter 8, we see that a pseudodifferential
operator is Heisenberg-smooth, because $\|a(x,D)\|$ depends just on a finite
number of seminorms of $a\in CB^\infty(\R^{2n},A)$. Besides, 
for $T=a(x,D)$, 
we have $\partial^\beta_z\partial^\gamma_\zeta T_{z,\zeta} =
\ope(\partial^\beta_z\partial^\gamma_\zeta a_{z,\zeta})$,  where 
$a_{z,\zeta}(x,\xi)= a(x+z,\xi+\zeta)$, \,\,\,$\beta,\gamma\in\N^n$ (for
proving these results, we just need to do some computations which
we can check in proposition 4.7 of \cite{ire}).        }\end{ob}

\begin{nota} {\rm Let $\HE$ be the subset of $\op$ 
formed by the 
Heisenberg-smooth operators. 
We have that $\ope:\CtwoRA\la\HE$ is a well
defined, injective application. For $A=\C$, we have that $\ope$ is a bijection, 
\cite{co-2}. }\end{nota}

\section{\bf Pseudodifferential operators that commute with $R_G$.}

Let us consider here the right regular representation of $\CBA$ for the
deformed product:

\[R_GF = F\times_JG.\]

\begin{lema}\label{aproxL}  If an operator $T\in\op$ commutes with $R_\varphi$
for all $\varphi\in\sa$,  there exists a sequence $F_k$ in $E$ such that
$F_k\times_J\varphi$ converges to $T(\varphi)$, for all $\varphi\in\sa$.   \end{lema}

\pf Let us find, first, a sequence, $e_k$, such that, for all $\varphi\in\sa$, 
$e_k\times_J\varphi \la
\varphi$ (convergence in the $\|\cdot\|_2$ norm).  
Since $A$ is separable, it has an approximate unit $(u_k)_{k\in\N}$.
For each $k\in\N$, let us consider a C$^\infty-$function $\phi_k:\R^n\rightarrow A$,
with support the set $\{x\in\R^n / \|x\|\leq\frac{1}{k}\}$ such that
$\int\phi_k(x)dx = u_k$.
Then, 
let $e_k = {\cal F}^{-1}(\phi_k)$, where ${\cal F}$ is the Fourier transform
on $\sa$ (for details, see  proposition 2.5 of \cite{ire}). 

Then,
since $R_\varphi$ is a continuous operator on $E$,  letting 
$F_k = Te_k \in E$,  $R_\varphi(F_k)$ is well defined  and
we have $F_k\times_J\varphi = R_\varphi Te_k = 
TR_\varphi e_k = T(e_k\times_J\varphi).$ Hence, since $e_k\times_J\varphi \longrightarrow 
\varphi$, we have $F_k\times_J\varphi \rightarrow T\varphi$,
for all $\varphi\in\sa$.\,\ai   



\begin{prop}\label{heissmooth}  If $T$ is an operator in $\op$ which
is such that $[T,R_\varphi]=0 \,\, \forall \varphi\in\sa$, then 
$T_{z,\zeta}=T_{z-J\zeta,0}$.\end{prop}

\pf  Since $R_\varphi$ is continuous, for any $F\in E$ and $\varphi\in\sa$,
we may write $L_{F}(\varphi) = R_\varphi(F) $ thus defining
 $L_{F}$  as an operator from $\sa$ to $E$. 

It is easy to see that
\[ E^{-1}_{z,\zeta}L_f E_{z,\zeta} = E^{-1}_{z-J\zeta,0}L_f E_{z-J\zeta,0},\] for 
$f\in\sa$. Then, using that $E_{z,\zeta}$ leaves $\sa$ invariant (here we are
writing $E_{z,\zeta}$ for $E_{z,\zeta,0}$, earlier defined), we get

\begin{equation}\label{EsyLs}
 E^{-1}_{z,\zeta}L_F E_{z,\zeta} = E_{z-J\zeta,0}L_F E_{z-J\zeta,0}, \end{equation} 
for $F\in E$,
so we have that $(L_F)_{z,\zeta}=(L_F)_{z-J\zeta,0}$. By lemma \ref{aproxL}, there is a sequence $F_k$ in $E$ such that, for all $\varphi\in\sa$,
$\dps\lim_{k\rr\infty}(L_{F_k})_{z,\zeta}(\varphi) = T_{z,\zeta}(\varphi)$ (by equation (\ref{EsyLs})),  so that
$T_{z,\zeta}=T_{z-J\zeta,0}$ .$\ai$

\begin{corol} \label{trans_smooth} If $T\in\op$ is such as in proposition \ref{heissmooth} and is translation-smooth,
then $T$ is Heisenberg-smooth. \end{corol}

\begin{lema}\label{operb}  Given $a\in\CtwoRA$, let $b = 
\dps\Pi_{j=1}^n (1+\partial_{ y_j})^2(1+\partial_{\xi_j})^2a$, and
$\gamma(x)=\dps\Pi_{j=1}^n f(x_j)$, with 

\[f(x_j)=\left\{\begin{array}{l}x_je^{-x_j} \,\, \text{if}\,\, x_j\geq 0 \\
    0 \hspace{1cm}\text{otherwise} \end{array}\right.\] Then we have
$a(x,\xi) = \dps\int\gamma(-z)\gamma(-\zeta)b(x+z.\xi+\zeta)dzd\zeta.$     \end{lema}

\noindent In the scalar case, $A=\C$, we can see the proof in 
\cite{co-1}, chapter 8, corollary 2.4.
The same argument is valid for the general case. 

\begin{teo}\label{TEOR} Let $A$ be a $C^*-$algebra for which the above defined 
application $\ope:\CtwoRA\la\HE$ is a bijection. Then, given an operator $T\in\op$,
translation-smooth, that commutes with $R_\varphi$ for all $\varphi\in\sa$,
there exists a function $F$ in $\CBA$ such that $T=L_F$.                   \end{teo}

\pf Since $\ope$ is a bijection, and by corollary \ref{trans_smooth},
there exists $a\in\CtwoRA$ such that $T=\ope(a)$.
As in lemma \ref{operb}, let $b = \dps\Pi_{j=1}^n(1+\partial_{x_j})^2(1+\partial_{\xi_j})^2a$
and $B = \ope(b)$. Note that 
$ B_{z,\zeta} = \Pi^n_{j=1}(1+\partial_{z_j})^2(1+\partial_{\zeta_j})^2T_{z,\zeta},$ see remark \ref{scalar}.


Since $TR_\varphi = R_\varphi T$, it is not difficult to see that $T_{z,\zeta}R\varphi =
R_\varphi T_{z,\zeta}$, for all $\varphi\in\sa$. Then, we have $[B,R_\varphi]=0$ 
for all $\varphi\in\sa$. So, by proposition \ref{heissmooth}, 
$B_{z,\zeta}=B_{z-J\zeta,0}$, so that $b(x+z,\xi+\zeta) = b(x+z-J\zeta,0)$.

By lemma \ref{operb}, we get $a(z,\zeta) = a(z-J\zeta,0)$. Choosing $F(z)=a(z,0)$, we
have $T=L_F$, with $F\in\CBA$, as was to be proved.$\ai$

\begin{ob}{\rm We have just proved that a pseudodifferential operator that commutes with all operators $G(x+JD), \,\,\,(G\in\sa)$, where $J$ is a fixed skew-symmetric matrix, 
is of the form $F(x-JD), \,\,\,F\in\sa$.                            } \end{ob}

\begin{ob}{\rm  If $A=\C$, since we have the Cordes' characterization \cite{co-2},
we can see that theorem \ref{TEOR} gives  us another proof of the main result of \cite{m-m},
without needing to apply the Schwartz kernel.           }\end{ob}

\bigskip

\centerline{\sc Acknowledgement}

This work was based on \cite{ire} and I benefited from many interchanges of
ideas, conversations, and suggestions about the redaction from my ex-advisor 
Prof. Severino T. Melo, to whom I remain deeply grateful.

\bigskip
{\sc Instituto de Matem\'atica e Estat\'\i stica, Universidade de S\~ao Paulo,
Caixa Postal 66281, 05315-970, S\~ao Paulo, Brazil.

{\it E-mail address:} {\bf marcela@ime.usp.br}

\end{document}